\documentclass[12pt,a4paper]{amsart}
\usepackage{amssymb,amsmath}

\headsep=1cm \numberwithin{equation}{section}
\newtheorem{theorem}{Theorem}[section]
\newtheorem{lemma}[theorem]{Lemma}
\newtheorem{proposition}[theorem]{Proposition}
\newtheorem{corollary}[theorem]{Corollary}

\theoremstyle{definition}

\theoremstyle{remark}

\newtheorem{example}[theorem]{Example}
\newtheorem{question}[theorem]{Question}

\newcommand{\Ass}{\operatorname{Ass}}

\newcommand{\grade}{\operatorname{grade}}

\newcommand{\pd}{\operatorname{pd}}
\newcommand{\injdim}{\operatorname{injdim}}
\newcommand{\Ext}{\operatorname{Ext}}
\newcommand{\Supp}{\operatorname{Supp}}

\newcommand{\Hom}{\operatorname{Hom}}
\newcommand{\Att}{\operatorname{Att}}

\newcommand{\depth}{\operatorname{depth}}

\newcommand{\Max}{\operatorname{Max}}

\newcommand{\BN}{\Bbb N}

\newcommand{\lo}{\longrightarrow}
\newcommand{\fm}{\frak{m}}
\newcommand{\fp}{\frak{p}}

\newcommand{\fa}{\frak{a}}

\begin{document}
\author[Divaani-Aazar, Sazeedeh and Tousi]{K. Divaani-Aazar, R. Sazeedeh and M. Tousi}
\title[Vanishing of generalized local cohomology]
{On vanishing of generalized local cohomology modules}

\address{K. Divaani-Aazar, Department of Mathematics, Az-Zahra University,
Vanak, Post Code 19834, Tehran, IRAN and Institute for Studies in
Theoretical Physics and Mathematics, P.O.Box 19395-5746, Tehran,
Iran} \email{kdivaani@ipm.ir}

\address{R. Sazeedeh, Department of Mathematics, Uromeiyeh University, Uromeiyeh,
IRAN and Institute of Mathematics, University for Teacher
Education, 599 Taleghani Avenue, Tehran 15614, IRAN.}

\address{M. Tousi, Institute for Studies in Theoretical Physics and
Mathematics, P. O. Box 19395-5746, Tehran, IRAN and Department of
Mathematics, Shahid Beheshti University, Tehran, IRAN.}
\email{mtousi@ipm.ir}



\thanks{This research was in part supported by a grant from IPM}

\begin{abstract}
Let $\fa$ denote an ideal of a $d$-dimensional Gorenstein local
ring $R$ and $M$ and $N$ two finitely generated $R$-modules with
$\pd M< \infty$. It is shown that $H^d_{\fa}(M,N)=0$ if and only
if $\dim \hat{R}\big/ \fa\hat{R}+\fp>0$ for all
$\fp\in\Ass_{\hat{R}}\hat{M}\cap\Supp_{\hat{R}}\hat{N}$.
\end{abstract}

\maketitle

\section{Introduction}
A generalization of local cohomology functors has been given by J.
Herzog in [{\bf 6}]. Let $\fa$ denote an ideal of a commutative
Noetherian ring $R$. For each $i\geq 0$, the functor
$H^i_{\fa}(.,.)$ defined by
$H^i_{\fa}(M,N)=\underset{n}{\varinjlim}\Ext_R^i(M/\fa^nM,N)$, for
all $R$-modules $M$ and $N$. Clearly, this notion is a
generalization of the usual local cohomology functor. The study of
this concept was continued in the articles [{\bf 8}], [{\bf
2}],[{\bf 9}], [{\bf 1}] and [{\bf 10}].

Two important type of theorems concerning local cohomology are
finiteness and vanishing results. We collect the known vanishing
results for generalized local cohomology in the following theorem.

\begin{theorem} Let $M$ and $N$ be two non-zero finitely
generated $R$-modules such that $\pd M<\infty$.

\begin{itemize}
\item[(i)] ([{\bf 9}, Theorem 3.7]) {\it Suppose $\dim N<\infty$.
Then $H^i_{\fa}(M,N)=0$ for all $i>\pd M+\dim (M\otimes_RN)$.}
\item[(ii)] ([{\bf 2}, Proposition 5.5]) {\it Let
$t=\grade_N(M/\fa M)=\inf\{i: \Ext_R^i(M/\fa M,N)\neq 0\}$. If
$t<\infty$, then $H^i_{\fa}(M,N)=0$ for all $i<t$ and
$H^t_{\fa}(M,N)\neq 0$.}
\item[(iii)] ([{\bf 9}, Theorem 2.5])
{\it $H^i_{\fa}(M,N)=0$ for all $i>ara(\fa)+\pd M$, where
$ara(\fa)$ the arithmetic rank of the ideal $\fa$ is the least
number of elements of $R$ required to generate an ideal which has
the same radical as $\fa$.} \item[(iv)] ([{\bf 8}, Theorem 2.3])
{\it Let $(R,\fm)$ be a local ring. Then $depth N$ is the least
integer $i$ such that $H^i_{\fm}(M,N)\neq 0$.}
\end{itemize}
\end{theorem}

As the main result of this paper, we generalize the
Lichtenbum-Hartshorne vanishing theorem to generalized local
cohomology in the certain case, where $R$ is Gorenstein. Namely,
we prove:

\begin{theorem} Let $\fa$ denote an ideal of a
$d$-dimensional Gorenstein local ring $(R,\fm)$. Let $M$ and $N$
be two finitely generated $R$-modules with $\pd M<\infty$. Then
$H^d_{\fa}(M,N)$ is an Artinian $R$-module. Moreover the following
are equivalent:
\begin{itemize}
\item[(i)] $H^d_{\fa}(M,N)=0.$ \item[(ii)]  {\it $\dim
\hat{R}\big/ \fa\hat{R}+\fp>0$ for all
$\fp\in\Ass_{\hat{R}}\hat{M}\cap\Supp_{\hat{R}}\hat{N}$.}
\end{itemize}
\end{theorem}

Having  1.1(i) in mind, one  may think that $\pd M+ \dim
(M\otimes_RN)$ or $\Max \{\pd M,\dim N\}$ is the last integer $i$
such that $H^i_{\fa}(M,N)\neq 0$. We provide examples which shows
that this is not true even in the case $R$ is local and $\fa$ is
the maximal ideal of $R$.

All rings considered in this paper are assumed to be commutative
Noetherian with identity. In our terminology we follow that of the
text book [{\bf 4}].

\section{Main result}
Let $\fa$ denote an ideal of a ring $R$. The generalized local
cohomology defined by
$$H^i_{\fa}(M,N)=\underset{n}{\varinjlim}\Ext_R^i(M/\fa^nM,N)$$
for all $R$-modules $M$ and $N$. Note that this is in fact a
generalization of the usual local cohomology, because if $M=R$,
then $H^i_{\fa}(R,N)=H^i_{\fa}(N)$. We use the following lemma
several times in this paper. Its proof is easy and we lift it to
the reader.

\begin{lemma} Let $M$ and $N$ be two $R$-modules. The
following are hold.
\begin{itemize}
\item[(i)] {\it Let $0\lo N\lo E^{\cdot}$ be an injective
resolution of $N$. Then $$H^i_{\fa}(M,N)\cong
H^i((\Gamma_{\fa}(\Hom_R(M,E^{\cdot})))\cong
H^i(H^0_{\fa}(M,E^{\cdot})).$$ Moreover, if $M$ is finitely
generated, then $H^i_{\fa}(M,N)\cong
H^i(\Hom_R(M,\Gamma_{\fa}(E^{\cdot})))$}.
\item[(ii)]  {\it If
$f:R\lo S$ is a flat ring homomorphism, then
$$H_{\fa}^i(M,N)\otimes_R S\cong H^i_{\fa S}(M\otimes_RS,
N\otimes_RS).$$}
\end{itemize}
\end{lemma}

\begin{theorem} Let $\fm$ be a maximal ideal of $R$ and $M,N$
be two finitely generated $R$-modules. Then $H^i_{\fm}(M,N)$ is
Artinian for all $i\geq 0$.
\end{theorem}

\vspace*{0.3cm} {\bf Proof.} Let $0\lo N\lo E^{\cdot}$ be a
minimal injective resolution of $N$. By 2.1 (i), it follows that
$H^i_{\fm}(M,N)=H^i(\Hom_R(M,\Gamma_{\fm}(E^{\cdot})))$. Denote
the $k$-th term of $E^{\cdot}$ by $E^k$. Because any subquotient
of an Artinian $R$-module is also Artinian, it is enough to show
that for each $k$, the module $\Hom_R(M,\Gamma_{\fm}(E^k))$ is
Artinian. One can see easily that for any prime ideal $\fp$ of
$R$,
$$\Hom_R(M,\Gamma_{\fm}(E(R/\fp)))=\begin{cases} 0 & ,\fp\neq \fm\\
\Hom_R(M,E(R/\fp))& ,\fp=\fm . \end{cases}$$ Thus
$\Hom_R(M,\Gamma_{\fm}(E^k))$ is equal to the direct sum of
$\mu^k(\fm,N)$ copies of $\Hom_R(M,E(R/\fm))$. Here
$\mu^k(\fm,N)=\dim_{R/\fm}(\Ext_R^k(R/\fm,N))_{\fm}$ is the $k$-th
Bass number of $N$ with respect to $\fm$,  which is clearly
finite. Next, it is easy to see that $\Hom_R(M,E(R/\fm))$ is
Artinian and so $\Hom_R(M,\Gamma_{\fm}(E^k))$ is Artinian as
required. $\Box$

The following is our technical tool throughout this paper.

\begin{proposition}  Let $\fa$ be an ideal
of a $d$-dimensional Gorenstein local ring $(R,\fm)$  and $M$ a
finitely generated $R$-module. Then the following statements hold.
\begin{itemize}
\item[(i)] {\it $H^i_{\fm}(M,R)=\begin{cases} \Hom_R(M,E(R/\fm)), & i=d \\
0 & , i\neq d \end{cases}.$ In particular, $H_m^d(R)=E(R/\fm).$}
\item[(ii)] {\it Assume that $\pd M<\infty$. For any $R$-module
$N$ and any $i>d$, $H^i_{\fa}(M,N)=0$ and so $H^d_{\fa}(M,.)$ is a
right exact functor.}
\end{itemize}
\end{proposition}

\vspace*{0.3cm} {\bf Proof.} (i) Let $0\lo R\lo E^{\cdot}$ be a
minimal injective resolution of $R$. Then
$E^i\cong\oplus_{htp=i}E(R/\fp)$ for each $0\leq i\leq d$ and
$E^i=0$ for all $i>d$, and so
$$\Hom_R(M,\Gamma_{\fm}(E^i))=\begin{cases} \Hom_R(M,E(R/\fm)), & i=d \\
0 & , i\neq d \end{cases}.$$ Hence
$H^i_{\fm}(M,R)=H^i(\Hom_R(M,\Gamma_{\fm}(E^{\cdot})))=\begin{cases}
\Hom_R(M,E(R/\fm)), & i=d \\ 0 &, i\neq d \end{cases}.$ For $M=R$,
we have $H^d_{\fm}(R)=H^d_{\fm}(R,R)=\Hom_R(R,E(R/\fm))=E(R/\fm)$.

(ii) Because $\injdim R=d$, it follows that $H^i_{\fa}(M,R)=0$ for
all $i>d$. Now, by decreasing induction on $i>d$, we show that
$H^i_{\fa}(M,N)=0$ for all finitely generated $R$-modules $N$.
This will complete the proof, because any $R$-module is the direct
limit of a direct system consisting of finitely generated
$R$-modules and the functor $H^i_{\fa}(M,\cdot)$, commutes with
direct limits. The claim clearly holds for $i=\pd M+\dim R+1$ by
1.1 (i). Assume that $i>d$ and that the claim holds for $i+1$.
Now, we prove it for $i$. Let $N$ be a finitely generated
$R$-module. There is a short exact sequence
$$0\lo K\lo F\lo N\lo 0,$$ where $F$ is a finitely generated free
$R$-module. We deduce the long exact sequence
$$\dots\lo H^i_{\fa}(M,K)\lo H^i_{\fa}(M,F)\lo H^i_{\fa}(M,N) \lo
H^{i+1}_{\fa}(M,K)\lo\dots .$$ By assumption
$H^{i+1}_{\fa}(M,K)=0$. On the other hand, we have
$H^i_{\fa}(M,F)=0$, because the functor $H^i_{\fa}(M,.)$ is
additive  and $H^i_{\fa}(M,R)=0$. Thus $H^i_{\fa}(M,N)=0$. $\Box$

The theory of attached prime ideals for Artinian modules is dual
of the theory of primary decomposition for Noetherian modules. For
an account of this theory, we refer the reader to [{\bf 4},
Chapter 7]. The following may be regarded as a generalization of
Grothendieck non-vanishing theorem, in the case $R$ is Gorenstein.

\begin{lemma} Let $(R,\fm)$ be $\fa$ $d$-dimensional
Gorenstein local ring and let $M$ and $N$ be two finitely
generated $R$-modules with $\pd M<\infty$. Then
$\Att_R(H^d_{\fm}(M,N))=\Ass_RM\cap\Supp_RN$. In particular,
$H^d_{\fm}(M,N)\neq 0$ if and only if
$\Ass_RM\cap\Supp_RN\neq\emptyset$.
\end{lemma}

\vspace*{0.3cm} {\bf Proof.} Since the functor $H^d_{\fm}(M,.)$ is
right exact, it follows from 2.3 that $H^d_{\fm}(M,N)\cong
N\otimes_R H^d_{\fm}(M,R)\cong N\otimes_R \Hom_R(M,E(R/\fm))$.
Thus $H^d_{\fm}(M,N)\cong \Hom_R(\Hom_R(N,M),E(R/\fm))$, by [{\bf
7}, Lemma 3.60]. For a finitely generated $R$-module $C$, it is
known and one can check easily that
$\Att_R(\Hom_R(C,E(R/\fm)))=\Ass_RC$. Therefore, by [{\bf 3},
p.267, Proposition 10],
$$\Att_R(H^d_m(M,N))=\Ass_R(\Hom_R(N,M))=\Supp_R N\cap\Ass_R M.$$
 The last assertion follows immediately, because
for an Artinian $R$-module $A$, $\Att_RA$ is empty if and only if
$A$ is zero. $\Box$

Let $\fa$ be an ideal of $R$. For an Artinian $R$-module $A$, we
put $<\fa>A=\cap_{n\in \BN}\fa^nA$.

\begin{theorem} Let the situation be as
in 2.4. Let $\fa$ be an ideal of $R$. Then there is a natural
isomorphism
$$H^d_{\fa}(M,N)\cong
H^d_{\fm}(M,N)/\sum_{n\in\BN}<\fm>(0:_{H^d_{\fm}(M,N)}\fa^n).$$ In
particular, $H^d_{\fa}(M,N)$ is Artinian.
\end{theorem}

\vspace*{0.3cm} {\bf Proof.} Denote $E(R/\fm)$, by $E$. It follows
from the local duality theorem [{\bf 4}, 11.2.5] that
$$\Ext_R^d(M/\fa^nM,R)\cong \Hom_R(H^0_{\fm}(M/\fa^nM),
E).$$ Let $A=\Hom_R(M,E)$. For a fixed integer $n$, let
$t(n)\in\BN$ be such that
$H^0_{\fm}(M/\fa^nM)=\Hom_R(R/\fm^{t(n)},M/\fa^nM)$ and
$<\fm>(0:_A\fa^n)=\fm^{t(n)}(0:_A\fa^n)$. Then by [{\bf 7}, Lemma
3.60], it follows that
$$\mathrm{Hom}_R(H^0_{\fm}(M/\fa^nM),E)\cong
R/\fm^{t(n)}\otimes_R\Hom_R(M/\fa^nM,E)$$
$$\cong R/\fm^{t(n)}\otimes_R\Hom_R(R/\fa^n,A)\cong
(0:_A\fa^n)/<\fm>(0:_A\fa^n).$$ One can check easily that
$\underset{n}{\varinjlim}(0:_A\fa^n)/<\fm>(0:_A\fa^n)\cong
A/\sum_{n\in\BN}<\fm>(0:_A\fa^n) .$ Hence $$H^d_{\fa}(M,R)\cong
A/\sum_{n\in\BN}<\fm>(0:_A\fa^n).$$ It follows by [{\bf 5}, Lemma
3.1], that $H^d_{\fa}(M,R)\otimes_RN\cong
(A\otimes_RN)/\sum_{n\in\BN}<\fm>(0:_{A\otimes_R N}\fa^n)$. But
2.3 (ii) implies that $H^d_{\fa}(M,R)\otimes_RN\cong
H^d_{\fa}(M,N)$ and $A\otimes_RN\cong H_{\fm}^d(M,N)$. Note that
$A\cong H^d_{\fm}(M,R)$, by 2.3 (i). This finishes the proof.
$\Box$

The following is an extension of the Lichetenbum-Hartshorne
vanishing theorem for generalized local cohomology.

\begin{corollary} Let $\fa$ denote an ideal of a
$d$-dimensional Gorenstein local ring $(R,\fm)$. Let $M$ and $N$
be two finitely generated $R$-modules with $\pd M<\infty$. Then
$\Att_{\hat{R}}(H^d_{\fa}(M,N))=\{\fp\in\Ass_{\hat{R}}\hat{M}\cap\Supp_{\hat{R}}\hat{N}:
\dim \hat{R}/\fa\hat{R}+\fp=0\}$. Thus the following statements
are equivalent:
\end{corollary}

\vspace*{0.3cm}
\begin{itemize}
\item[(i)] $H^d_{\fa}(M,N)=0.$
\item[(ii)] {\it $\dim
\hat{R}/\fa\hat{R}+\fp>0,$ for all
$\fp\in\Ass_{\hat{R}}\hat{M}\cap\Supp_{\hat{R}}\hat{N}$.}
\end{itemize}
\vspace*{0.3cm} {\bf Proof.} Let $A=H^d_{\fm}(M,N)$ and
$B=\sum_{n\in\BN}<\fm>(0:_{H^d_{\fm}(M,N)}\fa^n)$. Let
$A=\sum_{n\in\BN}A_i$ be a minimal secondary representation of the
Artinian $\hat{R}$-module $A$, where $A_i$ is $\fp_i$-secondary.
 We may assume that $\dim \hat{R}/\fa\hat{R}+\fp_i>0$ for $i=1,\dots ,k$
 and  that $\dim \hat{R}/\fa\hat{R}+\fp_i=0$ for $i=k+1,\dots ,t$.
 Then by [{\bf 5}, Theorem 2.8], $\sum_{i=1}^kA_i$ is a
minimal secondary representation of $B$. It is easy to see that
$A/B=\sum_{i=k+1}^t(A_i+B)/B$ is a minimal secondary
representation of $A/B$. Therefore, it follows from 2.4 and 2.5
that
$$\Att_{\hat{R}}(H^d_{\fa}(M,N))=\{\fp\in\Ass_{\hat{R}}\hat{M}\cap\Supp_{\hat{R}}\hat{N}:
\dim \hat{R}/\fa\hat{R}+\fp=0\}.$$ Note that because
$H^d_{\fm}(M,N)$ is an Artinian $R$-module, we have
$$H^d_{\fm}(M,N) \cong H^d_{\fm}(M,N)\otimes_R{\hat{R}} \cong
H^d_{\fm \hat{R}}(\hat{M}, \hat{N}).\Box $$

\vspace*{0.3cm}

\begin{question} Let $(R,\fm)$ be a local ring and $M$ and $N$
two finitely generated $R$-modules with $\pd M<\infty$. Describe
the last integer $i$ such  that $H^i_{\fm}(M,N)\neq 0$.
\end{question}

The following examples shows that the above mentioned integer is
neither  $\pd M + \dim(M \otimes_RN)$ nor $\Max \{\pd M,\dim N
\}$.

\begin{example} Suppose that$(R,\fm)$ is a regular local ring of
dimension $d$.
\begin{itemize}
\item[(i)] Suppose that $d>1$. Let $\fp\neq \fm$ be a non-zero
prime ideal and $x$ a non-zero element in $\fp$. Set $N=R/xR$ and
$M=R/\fp$. Then by 2.4,
$$\Att_R(H^d_{\fm}(M,N))=\Ass_RM\cap\Supp_RN=\{\fp\}.$$
Hence $H^d_{\fm}(M,N)\neq 0$. On the other hand, since by
Auslander-Buchsbaum formula $\pd M=\depth R-\depth M$, we have
$\Max\{\pd M,\dim N\}=d-1$. \item[(ii)] Suppose $M\neq 0$ is a
non-Cohen-Macaulay finitely generated $R$-module. Then
$$\depth R=\pd M+ \depth M<\pd M+ \dim(M \otimes_RR)=l.$$ Thus  $H^l_{\fm}(M,R)=
0$, by 2.3(i).
\end{itemize}
\end{example}


\end{document}